\documentstyle[10pt]{amsart}

\newtheorem{prop}{Proposition}

\newtheorem{th}{Theorem}
\newtheorem{lemma}{Lemma}
\newtheorem{cor}{Corollary}

\title{Almost complex structures that are harmonic maps}

\author{Johann Davidov, Absar Ul Haq, Oleg Mushkarov}
\thanks{The first and the third named authors are partially supported by  the National Science
Fund, Ministry of Education and Science of Bulgaria under contract
DFNI-I 02/14. }

\address{Johann Davidov\\Institute of Mathematics and Informatics \\
Bulgarian Academy of Sciences\\ Acad. G.Bonchev Str. Bl.8\\ 1113
Sofia\\ Bulgaria.} \email{jtd@@math.bas.bg}

\address{Absar ul Haq\\Department of Mathematics\\University of Management and Technology Lahore\\
Sialkot Campus\\ Pakistan.}\email{absar.ulhaq@@skt.umi.edu.pk}

\address{Oleg Mushkarov \\Institute of Mathematics and Informatics \\
Bulgarian Academy of Sciences\\ Acad. G.Bonchev Str. Bl.8\\ 1113
Sofia\\ Bulgaria {\it and} South-West University "Neofit Rilski"
\\2700 Blagoevgrad\\Bulgaria.}\email{muskarov@@math.bas.bg}

\begin{document}

\begin{abstract}
We find geometric conditions on a four-dimensional almost Hermitian
manifold under which the almost complex structure is a harmonic map
or a minimal isometric imbedding of the manifold into its twistor
space.

\vspace{0,1cm} \noindent 2010 {\it Mathematics Subject
Classification}. Primary 53C43, Secondary 58E20, 53C28

\vspace{0,1cm} \noindent {\it Key words: almost complex structures,
twistor spaces, harmonic maps}
\end{abstract}

\thispagestyle{empty}

\maketitle
\vspace{0.5cm}

\section{Introduction}

Recall that an almost complex structure on a Riemannian manifold
$(M,g)$ is called almost Hermitian if it is $g$-orthogonal. If a
Riemannian manifold admits an almost Hermitian structure $J$, it has
many such structures. One way to see this is to consider the twistor
bundle $\pi:{\cal Z}\to M$ whose fibre at a point $p\in M$ consists
of all $g$-orthogonal complex structures $I_p:T_pM\to T_pM$
($I_p^2=-Id$) on the tangent space of $M$ at $p$ yielding the same
orientation as $J_p$ . The fibre is the compact Hermitian symmetric
space $SO(2n)/U(n)$ and its standard metric
$-\frac{1}{2}Trace\,I_1\circ I_2$ is K\" ahler-Einstein. The twistor
space admits a natural Riemannian metric $h$ such that the
projection map $\pi:({\cal Z},h)\to (M,g)$ is a Riemannian
submersion with totally geodesic fibres. Consider $J$ as a section
of the bundle $\pi:{\cal Z}\to M$ and take a section $V$ with
compact support $K$ of the bundle $J^{\ast}{\cal V}\to M$, the
pull-back under $J$ of the vertical bundle ${\cal V}\to {\cal Z}$.
There exists $\varepsilon>0$ such that,  for every point $I$ of the
compact set $J(K)$,  the exponential map $exp_{I}$ is a
diffeomorphism of the $\varepsilon$-ball in $T_I{\cal Z}$. The
function $||V||_{h}$ is bounded on $M$, so there exists a number
$\varepsilon'>0$ such that $||tV(p)||_{h}<\varepsilon$ for every
$p\in M$ and $t\in(-\varepsilon',\varepsilon')$. Set
$J_t(p)=exp_{J(p)}[tV(p)]$ for $p\in M$ and
$t\in(-\varepsilon',\varepsilon')$. Then $J_t$ is a section of
${\cal Z}$, i.e. an almost Hermitian structure on $(M,g)$ (such that
$J_t=J$ on $M\setminus K$).

Thus it is natural to seek for "reasonable" criteria that
distinguish some of the almost Hermitian structures on a given
Riemannian manifold (cf., for example, \cite{CG,W1,W2,BLS}).
Motivated by the harmonic map theory, C. Wood \cite{W1,W2} has
suggested to consider as "optimal" those almost Hermitian structures
$J:(M,g)\to ({\cal Z},h)$ that are critical points of the energy
functional under variations through sections of ${\cal Z}$, i.e.
that are harmonic sections of the twistor bundle. In general, these
critical points are not harmonic maps, but, by analogy, they are
referred to as "harmonic almost complex structures"  in
\cite{W1,W2}. It  is more appropriate in the context of this article
to call  such structures "harmonic sections", a term used also in
\cite{W1}.

The almost Hermitian structures that are critical points of the
energy functional under variations through all maps $M\to{\cal Z}$
are genuine harmonic maps and the purpose of this paper is to find
geometric conditions on a four-dimensional almost Hermitian manifold
$(M,g,J)$ under which the almost complex structure $J$ is a harmonic
map of $(M,g)$ into $({\cal Z},h)$. We also find conditions for
minimality of the submanifold $J(M)$ of the twistor space. As   is
well-known, in dimension four, there are three basic classes in the
Gray-Hervella classification \cite{GH} of almost Hermitian
structures - Hermitian, almost K\"ahler (symplectic) and K\"ahler
structures. If $(g,J)$ is K\"ahler, the map $J:(M,g)\to ({\cal
Z},h)$ is a totally geodesic isometric imbedding.  In the case of a
Hermitian structure, we express the conditions for harmonicity and
minimality of $J$ in terms of the Lee form, the Ricci and star-Ricci
tensors of $(M,g,J)$, while for an almost K\"ahler structure the
conditions are in terms of the Ricci, star-Ricci and Nijenhuis
tensors. Several examples illustrating these results are discussed
in the last section of the paper, among them a Hermitian structure
that is a harmonic section of the twistor bundle and a minimal
isometric imbedding in it but not a harmonic map.

\smallskip

\noindent {\bf Acknowledgment}. The authors would like to thank the
referee for his/her remarks.

\section{Preliminaries}

Let $(M,g)$ be an oriented  Riemannian manifold of dimension four.
The metric $g$ induces a metric on the bundle of two-vectors
$\pi:\Lambda^2TM\to M$ by the formula
$$
g(v_1\wedge v_2,v_3\wedge v_4)=\frac{1}{2}det[g(v_i,v_j)].
$$
The Levi-Civita connection of $(M,g)$ determines a connection on the
bundle $\Lambda^2TM$, both denoted by $\nabla$, and the
corresponding curvatures are related by
$$
R(X\wedge Y)(Z\wedge T)=R(X,Y)Z\wedge T+Z\wedge R(X,Y)T
$$
for $X,Y,Z,T\in TM$. The  curvature operator ${\cal R}$ is the
self-adjoint endomorphism of $\Lambda ^{2}TM$ defined by
$$
 g({\cal R}(X\land Y),Z\land T)=g(R(X,Y)Z,T)
$$
Let us note that we adopt the following definition for the
curvature tensor $R$ :
$R(X,Y)=\nabla_{[X,Y]}-[\nabla_{X},\nabla_{Y}]$.

 The Hodge star operator defines an endomorphism $\ast$ of
$\Lambda^2TM$ with $\ast^2=Id$. Hence we have the orthogonal
decomposition
$$
\Lambda^2TM=\Lambda^2_{-}TM\oplus\Lambda^2_{+}TM
$$
where $\Lambda^2_{\pm}TM$ are the subbundles of $\Lambda^2TM$
corresponding to the $(\pm 1)$-eigenvalues of the operator $\ast$.

\smallskip

Let $(E_1,E_2,E_3,E_4)$ be a local oriented orthonormal frame of
$TM$. Set
\begin{equation}\label{s-basis}
s_1=E_1\wedge E_2 + E_3\wedge E_4, \quad s_2=E_1\wedge E_3+
E_4\wedge E_2, \quad s_3=E_1\wedge E_4+ E_2\wedge E_3.
\end{equation}
Then $(s_1,s_2,s_3)$ is a local orthonormal frame of
$\Lambda^2_{+}TM$ defining an orientation on $\Lambda^2_{+}TM$,
which does not depend on the choice of the frame
$(E_1,E_2,E_3,E_4)$.

\smallskip

For every $a\in\Lambda ^2TM$, define a skew-symmetric endomorphism
$K_a$ of $T_{\pi(a)}M$ by

\begin{equation}\label{cs}
g(K_{a}X,Y)=2g(a, X\wedge Y), \quad X,Y\in T_{\pi(a)}M.
\end{equation}
Note that, denoting by $G$ the standard metric
$-\frac{1}{2}Trace\,PQ$ on the space of skew-symmetric
endomorphisms, we have $G(K_a,K_b)=2g(a,b)$ for $a,b\in \Lambda
^2TM$. If $\sigma\in\Lambda^2_{+}TM$ is a unit vector, then
$K_{\sigma}$ is a complex structure on the vector space
$T_{\pi(\sigma)}M$ compatible with the metric and the orientation of
$M$. Conversely, the $2$-vector $\sigma$ dual to one half of the
fundamental $2$-form of such a complex structure is a unit vector in
$\Lambda^2_{+}TM$. Thus the unit sphere subbunlde ${\mathcal Z}$ of
$\Lambda^2_{+}TM$ parametrizes the complex structures on the tangent
spaces of $M$ compatible with its metric and orientation. This
subbundle is called the twistor space of $M$.

\smallskip

The Levi-Civita connection $\nabla$ of $M$ preserves the bundles
$\Lambda^2_{\pm}TM$, so it induces a metric connection on each of
them denoted again by $\nabla$. The  horizontal distribution of
$\Lambda^2_{+}TM$ with respect to $\nabla$ is tangent to the
twistor space ${\cal Z}$. Thus we have the decomposition $T{\cal
Z}={\cal H}\oplus {\cal V}$ of the tangent bundle of ${\cal Z}$
into horizontal and vertical components. The vertical space ${\cal
V}_{\tau}=\{V\in T_{\tau}{\cal Z}:~ \pi_{\ast}V=0\}$ at a point
$\tau\in{\cal Z}$ is the tangent space to the fibre of ${\cal Z}$
through $\tau$. Considering $T_{\tau}{\cal Z}$ as a subspace of
$T_{\tau}(\Lambda^2_{+}TM)$ (as we shall always do), ${\cal
V}_{\tau}$ is the orthogonal complement of $\tau$ in
$\Lambda^2_{+}T_{\pi(\tau)}M$. The map $V\ni{\cal V}_{\tau}\to
K_{V}$ gives an identification of the vertical space with the
space of skew-symmetric endomorphisms of $T_{\pi(\tau)}M$ that
anti-commute with $K_{\tau}$. Let $s$ be a local section of ${\cal
Z}$ such that $s(p)=\tau$ where $p=\pi(\tau)$. Considering $s$ as
a section of $\Lambda^2_{+}TM$, we have $\nabla_{X}s\in{\cal
V}_{\tau}$ for every $X\in T_pM$ since $s$ has a constant length.
Moreover, $X^h_{\tau}=s_{\ast}X-\nabla_{X}s$ is the horizontal
lift of $X$ at ${\tau}$.

\smallskip

Denote by $\times$ the usual vector cross product on the oriented
$3$-dimensional vector space $\Lambda^2_{+}T_pM$, $p\in M$,
endowed with the metric $g$. Then it is easy to check that
\begin{equation}\label{r-r}
g(R(a)b,c)=g({\cal R}(b\times c),a)
\end{equation}
for $a\in\Lambda^2T_pM$, $b,c\in\Lambda^2_{+}T_pM$. It is also easy to show that for every
$a,b\in\Lambda^2_{+}T_pM$
\begin{equation}\label{com}
K_a\circ K_b=-g(a,b)Id+ K_{a\times b}.
\end{equation}

For every $t>0$, define a Riemannian metric $h_t$ by
$$
h_t(X^h_{\sigma}+V,Y^h_{\sigma}+W)=g(X,Y)+tg(V,W)
$$
for $\sigma\in{\cal Z}_{+}$, $X,Y\in T_{\pi(\sigma)}M$, $V,W\in{\cal V}_{\sigma}$.
\smallskip

The twistor space ${\cal Z}$ admits two natural almost complex
structures that are compatible with the metrics $h_t$. One of them
has been introduced by Atiyah, Hitchin and Singer who have proved
that it is integrable if and only if the base manifold is
anti-self-dual \cite{AHS}. The other one, introduced by Eells and
Salamon, although never integrable, plays an important role in
harmonic map theory \cite{ES}.

\smallskip

 The action of $SO(4)$ on $\Lambda^2{\Bbb R}^4$ preserves the
decomposition $\Lambda^2{\Bbb R}^4=\Lambda^2_{+}{\Bbb R}^4\oplus
\Lambda^2_{-}{\Bbb R}^4$. Thus, considering $S^2$ as the unit
sphere in $\Lambda^2_{+}{\Bbb R}^4$, we have an action of the
group $SO(4)$ on $S^2$. Then, if $SO(M)$ denotes  the principal
bundle of the oriented orthonormal frames on $M$, the twistor
space ${\cal Z}$ is the associated bundle $SO(M)\times_{SO(4)}
S^2$. It follows from the Vilms theorem (see, for example,
\cite[Theorem 9.59]{Besse}) that the projection map $\pi:({\cal
Z},h_t)\to (M,g)$ is a Riemannian submersion with totally geodesic
fibres (this can also be proved by a direct computation).

\smallskip

Let $(N,x_1,...,x_4)$ be a local coordinate system of $M$ and let
$(E_1,...,E_4)$ be an oriented orthonormal frame of $TM$ on $N$. If
$(s_1,s_2,s_3)$ is the local frame of $\Lambda^2_{+}TM$ defined by
(\ref{s-basis}), then $\widetilde x_{a}=x_{a}\circ\pi$,
$y_j(\tau)=g(\tau, (s_j\circ\pi)(\tau))$, $1\leq a \leq 4$, $1\leq
j\leq 3$, are local coordinates of $\Lambda^2_{+}TM$ on
$\pi^{-1}(N)$.

   The horizontal lift $X^h$ on $\pi^{-1}(N)$ of a vector field
$$
X=\sum_{a=1}^4 X^{a}\frac{\partial}{\partial x_{a}}
$$
is given by
\begin{equation}\label{hl}
X^h=\sum_{a=1}^4 (X^{a}\circ\pi)\frac{\partial}{\partial \widetilde{x}_{a}}
-\sum_{j,k=1}^3y_j(g(\nabla_{X}s_j,s_k)\circ\pi)\frac{\partial}{\partial y_k}.
\end{equation}
Hence
\begin{equation}\label{Lie-1}
[X^h,Y^h]=[X,Y]^h+\sum_{j,k=1}^3y_j(g(R(X\wedge Y)s_j,s_k)\circ\pi)\frac{\partial}{\partial y_k}
\end{equation}
for every vector fields $X$ and $Y$ on $N$. Let $\tau\in{\cal Z}$.
Using the standard identification
$T_{\tau}(\Lambda^2_{+}T_{\pi(\tau)}M)\cong
\Lambda^2_{+}T_{\pi(\tau)}M$,  we obtain from (\ref{Lie-1}) the
well-known formula
\begin{equation}\label{Lie-2}
[X^h,Y^h]_{\tau}=[X,Y]^h_{\tau}+R_{p}(X\wedge Y)\tau, \quad p=\pi(\tau).
\end{equation}

Denote by $D$ the Levi-Civita connection of $({\cal Z},h_t)$. Then,
using the  Koszul formula for the Levi-Civita connection, identity
(\ref{Lie-2}), and the fact that the fibers of ${\mathcal Z}$ are
totally geodesic submanifolds,  it is easy to show  the following.

\begin{lemma}\label{LC} {\rm (\cite{DM})}
If $X,Y$ are vector fields on $M$ and $V$ is a vertical vector
field on ${\cal Z}$, then
\begin{equation}\label{D-hh}
(D_{X^h}Y^h)_{\tau}=(\nabla_{X}Y)^h_{\tau}+\frac{1}{2}R_{p}(X\wedge Y)\tau,
\end{equation}
\begin{equation}\label{D-vh}
(D_{V}X^h)_{\tau}={\cal H}(D_{X^h}V)_{\tau}=-\frac{t}{2}(R_{p}(\tau\times V)X)^h_{\tau}
\end{equation}
where $\tau\in{\cal Z}$, $p=\pi(\tau)$, and ${\cal H}$ means "the
horizontal component".
\end{lemma}

\bigskip

Let $(M,g,J)$ be an almost Hermitian manifold of dimension four.
Define a section ${\frak J}$ of $\Lambda^2TM$ by
$$
g({\frak J},X\wedge Y)=\frac{1}{2}g(JX,Y),\quad X,Y\in TM.
$$
Note that the section $2{\frak J}$ is dual to the fundamental
$2$-form of $(M,g,J)$. Consider $M$ with the orientation induced
by the almost complex structure $J$. Then ${\frak J}$ takes its
values in the twistor space ${\cal Z}$ of the Riemannian manifold
$(M,g)$.

 Let ${\cal J}_1$ and ${\cal J}_2$ be the Atiyah-Hitchin-Singer and Eells-Salamon  almost complex
structures on the twistor space ${\cal Z}$. It is well-known (and
easy to see) that the map ${\frak J}:(M,J)\to ({\cal Z},{\cal
J}_1)$ is holomorphic if and only if the almost complex structure
$J$ is integrable, while ${\frak J}:(M,J)\to ({\cal Z},{\cal
J}_2)$ is holomorphic if and only if $J$ is symplectic (i.e.
$(g,J)$ is an almost K\"ahler structure).

  In this note we are  going to find geometric conditions  under
which the map ${\frak J}: (M,g)\to ({\cal Z},h_t)$ that represents
$J$ is harmonic.

Let ${\frak J}^{-1}T{\cal Z}\to M$ be the pull-back of the bundle
$T{\cal Z}\to {\cal Z}$ under the map ${\frak J}:M\to{\cal Z}$. Then
we can consider the differential ${\frak J}_{\ast}:TM\to T{\cal Z}$
as a section of the bundle $Hom(TM,{\frak J}^{-1}T{\cal Z})\to M$.
Denote by $D^{({\frak J})}$ the connection on ${\frak J}^{-1}T{\cal
Z}$ induced by the Levi-Civita connection $D$ on $T{\cal Z}$. The
Levi Civita connection $\nabla$ on $TM$ and the connection
$D^{({\frak J})}$ on ${\frak J}^{-1}T{\cal Z}$ induce a connection
$\widetilde\nabla$ on the bundle $Hom(TM,{\frak J}^{-1}T{\cal Z})$.
The map ${\frak J}: (M,g)\to ({\cal Z},h_t)$ is harmonic iff
$$
Trace\,\widetilde\nabla{\frak J}_{\ast}=0.
$$
(cf., for example, \cite{EL}). Recall also that the map ${\frak J}:
(M,g)\to ({\cal Z},h_t)$ is totally geodesic iff
$\widetilde\nabla{\frak J}_{\ast}=0$.

\begin{prop}\label{covder-dif}
For every  $X,Y\in T_pM$, $p\in M$,
$$
\begin{array}{c}
\widetilde\nabla {\frak J}_{\ast}(X,Y)
=\displaystyle{\frac{1}{2}}[\nabla^{2}_{XY}{\frak J}
-g(\nabla^{2}_{XY}{\frak J},{\frak J}){\frak J}(p)
+ \nabla^{2}_{YX}{\frak J}-g(\nabla^{2}_{YX}{\frak J},{\frak J}){\frak J}(p)\\[8pt]
\hspace{4cm}-t(R_p({\frak J}(p)\times\nabla_X{\frak
J})Y)^h_{\sigma}-t(R_p({\frak J}(p)\times\nabla_Y{\frak
J})X)^h_{\sigma}]
\end{array}
$$
where $\nabla^{2}_{XY}{\frak J}=\nabla_X\nabla_Y{\frak
J}-\nabla_{\nabla_XY}{\frak J}$ is the second covariant derivative
of the section ${\frak J}$ of $\Lambda^2_{+}TM$.
\end{prop}

{\bf Proof}. Extend $X$ and $Y$ to vector fields in a neighbourhood
of the point $p$. Take an oriented orthonormal frame $E_1,...,E_4$
near $p$ such that $E_3=JE_2$, $E_4=JE_1$, so ${\frak J}=s_3$.
Define coordinates $(\tilde x_a,y_j)$ as above by means of this
frame and a coordinate system of $M$ at $p$. Set
$$
\begin{array}{c}
V_1=\displaystyle{(1-y_2^2)^{-1/2}(y_3\frac{\partial}{\partial y_1}-y_1\frac{\partial}{\partial y_3})},\\[8pt]
V_2=\displaystyle{(1-y_2^2)^{-1/2}(-y_1y_2\frac{\partial}{\partial y_1}+(1-y_2^2)\frac{\partial}{\partial y_2}-y_2y_3\frac{\partial}{\partial y_3})}.
\end{array}
$$
Then $V_1,V_2$ is a $g$-orthonormal frame of vertical vector fields
in a neighbourhood of the point $\sigma={\frak J}(p)$ such that
$V_1\circ{\frak J}=s_1$, $V_2\circ{\frak J}=s_2$. Note also that
$[V_1,V_2]_{\sigma}=0$. This and the Koszul formula imply
$(D_{V_k}V_l)_{\sigma}=0$ since $D_{V_k}V_l$ are vertical vector
fields, $k,l=1,2$. Thus $D_{W}V_l=0$, $l=1,2$, for every vertical
vector $W$ at $\sigma$. Considering ${\frak J}$ as a section of
$\Lambda^2_{+}TM$, we have
$$
{\frak J}_{\ast}\circ Y=Y^h\circ{\frak J}+\nabla_{Y}{\frak
J}=Y^h\circ{\frak J}+\sum_{k=1}^2g(\nabla_{Y}{\frak
J},s_k)(V_k\circ{\frak J}),
$$
hence
$$
\begin{array}{c}
D^{({\frak J})}_{X}({\frak J}_{\ast}\circ Y)=(D_{{\frak J}_{\ast}X}Y^h)\circ{\frak J}+\sum_{k=1}^2g(\nabla_{Y}{\frak J},s_k)(D_{{\frak J}_{\ast}X}V_k)\circ{\frak J}\\[8pt]
+\sum_{k=1}^2[g(\nabla_X\nabla_{Y}{\frak J},s_k)+g(\nabla_{Y}{\frak
J},\nabla_{X}s_k)](V_k\circ{\frak J})
\end{array}
$$
This, in view of Lemma~\ref{LC},  implies
$$
\begin{array}{c}
D^{({\frak J})}_{X_p}({\frak J}_{\ast}\circ Y)=(\nabla_{X}Y)^h_{\sigma}+\displaystyle{\frac{1}{2}}R(X,Y)\sigma\\[8pt]
-\displaystyle{\frac{t}{2}(R_p({\frak J}(p)\times\nabla_X{\frak J})Y)^h_{\sigma}-\frac{t}{2}(R_p({\frak J}(p)\times\nabla_Y{\frak J})X)^h_{\sigma}}\\[8pt]
+\sum_{k=1}^2g(\nabla_X\nabla_{Y}{\frak
J},s_k)s_k(p)+\sum_{k=1}^2[g(\nabla_Y{\frak
J},s_k)[X^h,V_k]_{\sigma}+g(\nabla_{Y}{\frak
J},\nabla_{X}s_k)s_k(p)].
\end{array}
$$
An easy computation using (\ref{hl}) gives
$$
[X^h,V_1]_{\sigma}=g(\nabla_{X}s_1,s_2)s_2(p),\quad
[X^h,V_2]_{\sigma}=g(\nabla_{X}s_2,s_1)s_1(p).
$$
These identities imply
$$
\sum_{k=1}^2[g(\nabla_Y{\frak
J},s_k)[X^h,V_k]_{\sigma}+g(\nabla_{Y}{\frak
J},\nabla_{X}s_k)s_k(p)]=0.
$$
Since $g(\nabla_{X}{\frak J},{\frak J})=0$ for every $X\in T_pM$, we
have
$$
g(\nabla_{Y}\nabla_{X}{\frak J},{\frak J})=-g(\nabla_{X}{\frak
J},\nabla_{Y}{\frak J})=g(\nabla_{X}\nabla_{Y}{\frak J},{\frak J})
$$
Hence
$$
\sum_{k=1}^2g(\nabla_X\nabla_{Y}{\frak
J},s_k)s_k(p)=\nabla_{X}\nabla_{Y}{\frak
J}-\frac{1}{2}g(\nabla_{X}\nabla_{Y}{\frak J}
+\nabla_{Y}\nabla_{X}{\frak J},{\frak J}){\frak J}(p).
$$
It follows that
$$
\begin{array}{c}
\widetilde\nabla{\frak J}_{\ast}(X,Y)=D^{({\frak J})}_{X_p}({\frak J}_{\ast}\circ Y)
-(\nabla_{X}Y)^h_{\sigma}-\nabla_{\nabla_XY}{\frak J}\\[8pt]
=\displaystyle{\frac{1}{2}}[\nabla_{X}\nabla_{Y}{\frak J}-\nabla_{\nabla_XY}{\frak J}
-g(\nabla_{X}\nabla_{Y}{\frak J},{\frak J}){\frak J}(p)\\[8pt]
+\nabla_{Y}\nabla_{X}{\frak J}-\nabla_{\nabla_YX}{\frak J}-g(\nabla_{Y}\nabla_{X}{\frak J},{\frak J}){\frak J}(p)\\[8pt]
 -t(R_p({\frak J}(p)\times\nabla_X{\frak J})Y)^h_{\sigma}-t(R_p({\frak J}(p)\times\nabla_Y{\frak J})X)^h_{\sigma}].
\end{array}
$$

\begin{cor}
If $(M,g,J)$ is a K\"ahler surface, the map ${\frak J}:(M,g)\to
({\cal Z},h_t)$ is a totally geodesic isometric imbedding.
\end{cor}

\smallskip

\noindent {\bf Remark 1}.
 By a result of C. Wood \cite{W1,W2},  ${\frak J}$ is a harmonic section of the twistor space $({\cal
Z},h_t)\to (M,g)$ if and only if $[J,\nabla^{\ast}\nabla J]=0$ where
$\nabla^{\ast}\nabla$ is the rough Laplacian. Taking into account
that $\nabla^{\ast}\nabla J=-Trace\,\nabla^2 J$, one can see that
the latter condition is equivalent to
$$
g(Trace\,\nabla^2{\frak J},X\wedge Y-JX\wedge JY)=0,\quad X,Y\in TM,
$$
which is equivalent to
$
{\cal V} Trace\,\nabla^2{\frak J}=0.
$
Thus, by Proposition~\ref{covder-dif}, ${\frak J}$ is a harmonic
section if and only if
$$
{\cal V} Trace\,\widetilde\nabla{\frak J}_{\ast}=0.
$$

\section{Harmonicity of ${\frak J}$}

Let $\Omega(X,Y)=g(JX,Y)$ be the fundamental $2$-form of the
almost Hermitian manifold $(M,g,J)$. Denote by $N$ the Nijenhuis
tensor of $J$
$$N(Y,Z)=-[Y,Z]+[JY,JZ]-J[Y,JZ]-J[JY,Z].$$ It is
well-known (and easy to check) that
\begin{equation}\label{nJ}
2g((\nabla_XJ)(Y),Z)=d\Omega(X,Y,Z)-d\Omega(X,JY,JZ)+g(N(Y,Z),JX),
\end{equation}
for all $X,Y,Z\in TM$.

\subsection{The case of integrable $J$}

Suppose that the almost complex structure $J$ is integrable. This
is equivalent to $(\nabla_{X}J)(Y)=(\nabla_{JX}J)(JY)$, $X,Y\in
TM$ \cite[Corollary 4.2]{G}. Let $B$ be the vector field on $M$
dual to the Lee form $\theta=-\delta\Omega\circ J$ with respect to
the metric $g$. Then (\ref{nJ}) and the identity
$d\Omega=\theta\wedge\Omega$ imply the following well-known
formula
\begin{equation}\label{nablaJ}
2(\nabla_XJ)(Y)=g(JX,Y)B-g(B,Y)JX+g(X,Y)JB-g(JB,Y)X.
\end{equation}
We have
$$
g(\nabla_X{\frak J},Y\wedge
Z)=\displaystyle{\frac{1}{2}}g((\nabla_XJ)(Y),Z)
$$
and it follows that
\begin{equation}\label{nalpha}
\nabla_X{\frak J}=\displaystyle{\frac{1}{2}}(JX\wedge B+X\wedge JB).
\end{equation}
The latter identity implies
\begin{equation}\label{na2alpha}
\nabla^2_{XY}{\frak
J}=\displaystyle{\frac{1}{2}}[(\nabla_XJ)(Y)\wedge
B+Y\wedge(\nabla_XJ)(B)+JY\wedge\nabla_XB+Y\wedge J\nabla_XB].
\end{equation}

Now recall that  the $\ast$-Ricci tensor $\rho^{\ast}$ of the
almost Hermitian manifold $(M,g,J)$ is defined by
$$\rho^{\ast}(X,Y)=trace\{Z\to R(JZ,X)JY\}.$$ Note that
\begin{equation}\label{rho-star}
\rho^{\ast}(JX,JY)=\rho^{\ast}(Y,X),
\end{equation}
in particular $\rho^{\ast}(X,JX)=0$.

Denote by $\rho$ the Ricci tensor of the Riemannian manifold
$(M,g)$.

\begin{th}\label{harm-int}
Suppose that the almost complex structure $J$ is integrable. Then
the map ${\frak J}:(M,g)\to ({\cal Z},h_t)$ is harmonic if and only
if $d\theta$ is a $(1,1)$-form and $\rho(X,B)=\rho^{\ast}(X,B)$ for
every $X\in TM$.
\end{th}

{\bf Proof}. According to Proposition~\ref{covder-dif} the map
${\frak J}$ is harmonic if and only if
$$
{\cal V}Trace\,\widetilde\nabla{\frak J}_{\ast}={\cal
V}Trace\,\nabla^2{\frak J}=0
$$
and
$$
\pi_{\ast}{\cal H}Trace\,\widetilde\nabla{\frak
J}_{\ast}=Trace\,\{TM\ni X\to R({\frak J}\times\nabla_X{\frak
J})X\}=0.
$$

It follows from identity (\ref{na2alpha})  that for every $X,Y\in TM$
\begin{equation}\label{Tr}
\begin{array}{c}
4g(Trace\,\nabla^2{\frak J},X\wedge Y)=\\[8pt]
-g(\nabla_{JX}B,Y)+g(\nabla_{JY}B,X)-g(\nabla_{X}B,JY)+g(\nabla_{Y}B,JX)+||B||^2g(X,JY)\\[8pt]
=-d\theta(JX,Y)-d\theta(X,JY)+||B||^2g(X,JY).
\end{array}
\end{equation}
Take an orthonormal frame $E_1,E_2=JE_1,E_3,E_4=JE_3$ so that
${\frak J}=E_1\wedge JE_1+E_3\wedge JE_3$. Then (\ref{Tr}) implies
$$
4g(Trace\,\nabla^2{\frak J},{\frak J})g({\frak J},X\wedge
Y)=-||B||^2g(JX,Y).
$$
Therefore
\begin{equation}\label{vtr}
\begin{array}{c}
4g({\cal V}\,Trace\,\nabla^2{\frak J},X\wedge
Y)=4[g(Trace\,\nabla^2{\frak J},X\wedge Y)-
g(Trace\,\nabla^2{\frak J},{\frak J})g({\frak J},X\wedge Y)]\\[8pt]
=-d\theta(JX,Y)-d\theta(X,JY).
\end{array}
\end{equation}
Thus ${\cal V}\,Trace\,\nabla^2{\frak J}=0$ if and only if
$d\theta(JX,Y)+d\theta(X,JY)=0$, which is equivalent to $d\theta$
being of type $(1,1)$.

Identity (\ref{nalpha}) implies that
$$
\begin{array}{c}
{\frak J}\times\nabla_{X}{\frak J}=g(\nabla_{X}{\frak J},s_2)s_3-g(\nabla_{X}{\frak J},s_3)s_2=g(\nabla_{JX}{\frak J},s_3)s_3+g(\nabla_{JX}{\frak J},s_2)s_2\\[6pt]
 =\nabla_{JX}{\frak J}=-\displaystyle{\frac{1}{2}}(X\wedge B-JX\wedge JB).
\end{array}
$$
It follows that for every $X\in TM$
\begin{equation}\label{htr}
2g(\pi_{\ast}{\cal H}Trace\,\widetilde\nabla{\frak
J}_{\ast},X)=2\sum_{i=1}^4g(R({\frak J}\times\nabla_{E_i}{\frak
J})E_i,X)=-\rho(X,B)+\rho^{\ast}(X,B)
\end{equation}
This proves the theorem.

\smallskip

Remark 1 and  formula (\ref{vtr}) in the proof of
Theorem~\ref{harm-int} give the following.
\begin{cor}\label{Wood}
The map ${\frak J}:(M,g)\to ({\cal Z},h_t)$ representing an
integrable almost Hermitian structure $J$ on $(M,g)$ is a harmonic
section if and only if the $2$-form $d\theta$ is of type $(1,1)$.
\end{cor}


\noindent {\bf Remark 2}. Note that the $2$-form $d\theta$ of a
Hermitian surface $(M,g,J)$ is of type $(1,1)$ if and only if the
$\star$-Ricci tensor $\rho^{\ast}$ is symmetric.

Indeed, let $s=Trace\,\rho$ and $s^{\ast}=Trace\,\rho^{\ast}$ be
the scalar and $\ast$-scalar curvatures. Set
$$
L(X,Y)=(\nabla_{X}\theta)(Y) +\frac{1}{2}\theta(X)\theta(Y),\quad
X,Y\in TM.
$$
Formula (3.4) in \cite{V82} implies the following identity
(\cite{M01})
$$
\rho(X,Y)-\rho^{\ast}(X,Y)=\frac{1}{2}[L(JX,JY)-L(X,Y)]+\frac{s-s^{\ast}}{4}g(X,Y).
$$
It follows that, when $\rho$ is $J$-invariant, $\rho^{\ast}$ is
symmetric  if and only if $L(JX,JY)-L(JY,JX)=L(X,Y)-L(Y,X)$. This
identity is equivalent to $d\theta(JX,JY)=d\theta(X,Y)$, which means
that $d\theta$ is of type $(1,1)$.

\smallskip

\subsection{The case of symplectic $J$}

Recall that an almost Hermitian manifold is called almost K\"ahler
(or symplectic) if its fundamental $2$-form is closed.

Denote by $\Lambda^2_{0}TM$ the subbundle of $\Lambda^2_{+}TM$
orthogonal to ${\frak J}$ (thus $\Lambda^2_{0}T_pM={\cal V}_{{\frak
J}(p)}$).  Under this notation we have the following.

\begin{th}\label{harm-sympl}
Let $(M,g,J)$ be an almost K\"ahler $4$-manifold. Then the map
${\frak J}:(M,g)\to ({\cal Z},h_t)$ is harmonic if and only if the
$\ast$-Ricci tensor $\rho^{\ast}$ is symmetric and
$$
Trace\,\{\Lambda^2_{0}TM\ni\tau\to R(\tau)(N(\tau))\}=0.
$$
\end{th}

{\bf Proof}. The $2$-form $\Omega$ is harmonic since $d\Omega=0$ and
$\ast\Omega=\Omega$, so by the Weitzenb\"ock formula
$$
(Trace\,\nabla^2\Omega)(X,Y)=Trace\{Z\to (R(Z,Y)\Omega)(Z,X)-(R(Z,X)\Omega)(Z,Y)\},
$$
$X,Y\in TM$ (see, for example, \cite{EL}). We have
$$
\begin{array}{c}
(R(Z,Y)\Omega)(Z,X)=-\Omega(R(Z,Y)Z,X)-\Omega(Z,R(Z,Y)X)\\[8pt]
=g(R(Z,Y)Z,JX)+g(R(Z,Y)X,JZ).
\end{array}
$$
Hence
$$
(Trace\,\nabla^2\Omega)(X,Y)=\rho(Y,JX)-\rho(X,JY)+2\rho^{\ast}(X,JY).
$$
Thus
\begin{equation}\label{tr-sympl}
2g(Trace\,\nabla^2{\frak J},X\wedge
Y)=(Trace\,\nabla^2\Omega)(X,Y)=\rho(Y,JX)-\rho(X,JY)+2\rho^{\ast}(X,JY).
\end{equation}
Let $E_1,...,E_4$ be an orthonormal basis of a tangent space $T_pM$
such that $E_2=JE_1$, $E_4=JE_3$. Define $s_1,s_2,s_3$ by
(\ref{s-basis}). Then ${\frak J}(p)=s_1$ and ${\cal V}_{{\frak
J}(p)}=span\{s_2,s_3\}$. By Proposition~\ref{covder-dif} and
identity (\ref{tr-sympl})
\begin{equation}\label{vtr-sympl}
\begin{array}{c}
g(Trace\,\widetilde\nabla^2{\frak
J}_{\ast},s_2)=g(Trace\,\nabla^2{\frak J},s_2)=\rho^{\ast}(E_1,E_4)
-\rho^{\ast}(E_4,E_1),\\[6pt]
g(Trace\,\widetilde\nabla^2{\frak
J}_{\ast},s_3)=g(Trace\,\nabla^2{\frak
J},s_3)=-\rho^{\ast}(E_1,E_3)-\rho^{\ast}(E_2,E_4).
\end{array}
\end{equation}
It follows, in view of (\ref{rho-star}), that ${\cal
V}Trace\,\widetilde\nabla^2{\frak J}_{\ast}=0$ if and only if
$\rho^{\ast}(E_i,E_j)=\rho^{\ast}(E_j,E_i)$, $i,j=1,...,4$.

 In order to compute $\pi_{\ast}{\cal H}Trace\widetilde\nabla
{\frak J}_{\ast}=Trace\,\{TM\ni X\to R({\frak J}\times\nabla_X{\frak
J})X\}$ we first note that, by (\ref{nJ}),
$$
g((\nabla_XJ)(Y),Z)=\frac{1}{2}g(N(Y,Z),JX).
$$
Then
$$
g(\nabla_X{\frak J},Y\wedge
Z)=\frac{1}{2}g((\nabla_XJ)(Y),Z)=\frac{1}{4}g(N(Y,Z),JX)
$$
The Nijenhuis tensor $N(Y,Z)$ is skew-symmetric, so it induces a
linear map $\Lambda^2TM\to TM$ which we denote again by $N$. It
follows that, for every $a\in\Lambda^2TM$ and $X\in T_{\pi(a)}M$,
\begin{equation}\label{nalpha-sympl}
g(\nabla_X{\frak J}, a)=\frac{1}{4}g(N(a),JX).
\end{equation}

Take an orthonormal basis $E_1,...,E_4$ of $T_pM$, $p\in M$, with
$E_2=JE_1$, $E_4=JE_3$. Define $s_1$, $s_2$, $s_3$ by
(\ref{s-basis}). Then ${\frak J}=s_1$ and
$$
\begin{array}{c}
{\frak J}\times\nabla_{X}{\frak J}=g(\nabla_{X}{\frak J},s_2)s_3-g(\nabla_{X}{\frak J},s_3)s_2=-\displaystyle{\frac{1}{4}}[g(JN(s_2),X)s_3-g(JN(s_3),X)s_2]\\[6pt]
=\displaystyle{\frac{1}{4}}[g(N(s_3),X)s_3+g(N(s_2),X)s_2]=-\nabla_{JX}{\frak
J},
\end{array}
$$
in view of (\ref{nalpha-sympl})  and  the identities
\begin{equation}\label{prop-N}
N(JX,Y)=N(X,JY)=-JN(X,Y),\quad X,Y\in T_pM. \end{equation}
Therefore
$$
\begin{array}{c} g(Trace\,\{TM\ni X\to R({\frak J}\times\nabla_X{\frak J})X\},Y)=-\sum_{i=1}^4 g(R(\nabla_{JE_i}{\frak J})E_i,Y)\\[6pt]
=\sum_{i=1}^4[g(\nabla_{JE_i}{\frak J},s_2)g(R(s_2)Y,E_i)+g(\nabla_{JE_i}{\frak J},s_3)g(R(s_3)Y,E_i)]\\[6pt]
=g(\nabla_{JR(s_2)Y}{\frak J},s_2)+ g(\nabla_{JR(s_3)Y}{\frak J},s_3)\\[6pt]
=-\displaystyle{\frac{1}{4}}[g(N(s_2),R(s_2)Y)+g(N(s_3),R(s_3)Y)]\\[6pt]
=\displaystyle{\frac{1}{4}}[g(R(s_2)(N(s_2)),Y)+g(R(s_3)(N(s_3)),Y)].
\end{array}
$$
Thus
\begin{equation}\label{nhtr-sympl}
4g(Trace\,\{TM\ni X\to R({\frak J}\times\nabla_X{\frak
J})X\},Y)=g(Trace\,\{\Lambda^2_{0}TM\ni\tau\to
R(\tau)(N(\tau))\},Y).
\end{equation}
This proves the theorem.
\smallskip

\section {Minimality of ${\frak J}$}

The map ${\frak J}:M\to {\cal Z}$ is an imbedding and, in this
section, we discuss the problem when ${\frak J}(M)$ is a minimal
submanifold of $({\cal Z},h_t)$.

Let $D'$ be the Levi-Civita connection of the metric on ${\frak
J}(M)$ induced
 by the metric $h_t$ on ${\cal Z}$. Let $\Pi$ be the second fundamental form
 of the submanifold ${\frak J}(M)$. Then, as is well-known (and easy to see), for every vector fields $X$ and $Y$ on $M$
$$
\widetilde\nabla{\frak J}_{\ast}(X,Y)=D'_{{\frak J}_{\ast}X}({\frak
J}_{\ast}\circ Y\circ{\frak J}^{-1})+\Pi({\frak J}_{\ast}X,{\frak
J}_{\ast}Y)-{\frak J}_{\ast}(\nabla_{X}Y).
$$
Thus $\Pi({\frak J}_{\ast}X,{\frak J}_{\ast}Y)$ is the normal
component of $\widetilde\nabla{\frak J}_{\ast}(X,Y)$, in particular
${\frak J}(M)$ is a minimal submanifold if and only if the normal
component of $Trace\,\widetilde\nabla{\frak J}_{\ast}$ vanishes.

\subsection{The case of integrable $J$}

\begin{th}\label{min-int}
Suppose that the almost complex structure $J$ is integrable. Then
the map ${\frak J}:M\to ({\cal Z},h_t)$ is a minimal isometric
imbedding if and only if  $d\theta$ is a $(1,1)$ form and
$\rho(X,B)=\rho^{\ast}(X,B)$ for every $X\perp \{B,JB\}$.
\end{th}

{\bf Proof}. Let $p\in M$ and suppose that $B_p=0$. Then
$\nabla{\frak J}|_p=0$ by (\ref{nalpha}). Hence ${\frak
J}_{\ast}(X)=X^h_{{\frak J}(p)}$ for every $X\in T_pM$. Thus the
tangent space of ${\frak J}(M)$ at the point ${\frak J}(p)$ is the
horizontal space ${\cal H}_{{\frak J}(p)}$, while the normal space
is the vertical space ${\cal V}_{{\frak J}(p)}$. Let $E_1$,
$E_2=JE_1$, $E_3$, $E_4=JE_3$ be an orthonormal basis of $T_pM$ and
define $s_1,s_2,s_3$ by formula (\ref{s-basis}). Then ${\frak
J}(p)=s_1$ and ${\cal V}_{{\frak J}(p)}=span\{s_2,s_3\}$ Hence the
normal component at ${\frak J}(p)$ of $Trace\,\widetilde\nabla{\frak
J}_{\ast}$ vanishes if and only if
$$
g({\cal V}Trace\,\widetilde\nabla{\frak J}_{\ast},s_2)=g({\cal
V}Trace\,\widetilde\nabla{\frak J}_{\ast},s_3)=0.
$$
Applying (\ref{vtr}), we see that this is equivalent to
$$
d\theta(E_2,E_3)=-d\theta(E_1,E_4),\quad
d\theta(E_2,E_4)=d\theta(E_1,E_3).
$$
The latter identities are equivalent to $(d\theta)_p$ being of type
$(1,1)$.

Now assume that $B_p\neq 0$. Then we can find an orthonormal basis
of $T_pM$ of the form $E, JE, ||B_p||^{-1}B_p, ||B_p||^{-1}JB_p$. It
follows from (\ref{nalpha}) that
$$
E^h_{{\frak J}(p)}-\frac{4}{t||B_p||^2}\nabla_{E}{\frak J},\quad
(JE)^h_{{\frak J}(p)}-\frac{4}{t||B_p||^2}\nabla_{JE}{\frak J}
$$
is a $h_t$-orthogonal basis of the normal space of ${\frak J}(M)$ at
${\frak J}(p)$. Therefore, according to
Proposition~\ref{covder-dif}, the normal component of the vertical
part of $Trace\,\widetilde\nabla{\frak J}_{\ast}$ at $s(p)$ vanishes
if and only if
$$
g({\cal V}\,Trace\,\nabla^2{\frak J},\nabla_{E}{\frak J})=g({\cal
V}\,Trace\,\nabla^2{\frak J},\nabla_{JE}{\frak J})=0.
$$
It follows from (\ref{nalpha}) and (\ref{vtr}) that the latter
identities hold if and only if $d\theta(X,B)=d\theta(JX,JB)$ for
every $X\perp \{B_p,JB_p\}$, which is equivalent to $d\theta$ being
of type $(1,1)$. Moreover, by (\ref{htr}), the normal component of
the horizontal part of $Trace\,\widetilde\nabla{\frak J}_{\ast}$
vanishes if and only if
$$
-\rho(E,B)+\rho^{\ast}(E,B)=-\rho(JE,B)+\rho^{\ast}(JE,JB)=0,
$$
or, equivalently, $\rho(X,B)=\rho^{\ast}(X,B)$ for every $X\in T_pM$, $X\perp \{B_p,JB_p\}$.

\begin{cor}
If $J$ is integrable and the Ricci tensor is $J$-invariant then, the
map ${\frak J}:M\to ({\cal Z},h_t)$ is a minimal isometric
imbedding.
\end{cor}

{\bf Proof}. According to \cite[Lemma 1]{M01}, the Ricci tensor
$\rho$ is $J$-invariant if and only if
$$
\rho-\rho^{\ast}=\frac{s-s^{\ast}}{4}g.
$$
Thus $\rho^{\ast}$ is symmetric if $\rho$ is $J$-invariant, hence
$d\theta$ is of type $(1,1)$ by Remark 2. Moreover, clearly
$\rho(X,B)=\rho^{\ast}(X,B)$ for $X\perp B$. Thus the result
follows from Theorem~\ref{min-int}.

This proof and Corollary~\ref{Wood} give the following

\begin{cor}
If $J$ is integrable and the Ricci tensor is $J$-invariant, then
${\frak J}:(M,g)\to ({\cal Z},h_t)$  is a harmonic section.
\end{cor}

\noindent {\bf Remark 3}. By \cite[Theorem 2]{AG}, every compact
Hermitian surface with $J$-invariant Ricci tensor is locally
conformally K\"ahler, $d\theta=0$. Moreover, if its first Betti
number is even, it is globally conformally K\"ahler (\cite{V80}). It
is still unknown whether there are compact complex surfaces with
$J$-invariant Ricci tensor and odd first Betti number.

\subsection{The case of symplectic $J$}

Set ${\cal N}_p=span\{N(X,Y):~X,Y\in T_pM\}$, $p\in M$, so ${\cal
N}_p=N(\Lambda^2T_pM)$. Identity (\ref{prop-N}) implies that
$N(\Lambda^2_{-}T_pM)=0$ and $N({\frak J})=0$. Hence ${\cal
N}_p=N(\Lambda^2_{0}T_pM)$ is a $J$-invariant subspace of $T_pM$
of dimension $0$ or $2$.

\begin{th}\label{min-sympl}
Let $(M,g,J)$ be an almost K\"ahler $4$-manifold. Then the map
${\frak J}:M\to ({\cal Z},h_t)$ is a minimal isometric imbedding if
and only if the $\star$-Ricci tensor $\rho^{\ast}$ is symmetric and
for every $p\in M$
$$
Trace\,\{\Lambda^2_{0}T_pM\ni\tau\to R_p(\tau)(N(\tau))\}\in{\cal
N}_p.
$$
\end{th}

{\bf Proof}. Suppose that $N_p=0$ for a point $p\in M$. Then
$\nabla{\frak J}|_p=0$ by (\ref{nalpha-sympl}). Hence ${\frak
J}_{\ast}(X)=X^h_{{\frak J}(p)}$ for every $X\in T_pM$. Thus the
normal space of ${\frak J}(M)$ at the point ${\frak J}(p)$ is the
vertical space ${\cal V}_{{\frak J}(p)}$. Therefore the normal
component at ${\frak J}(p)$ of $Trace\,\widetilde\nabla{\frak
J}_{\ast}$ vanishes if and only if
$$
g(Trace\,\widetilde\nabla{\frak
J}_{\ast},s_2)=g(Trace\,\widetilde\nabla{\frak J}_{\ast},s_3)=0.
$$
where $s_2$, $s_3$ are defined via (\ref{s-basis}) by means of an
orthonormal basis $E_1,...,E_4$ of $T_pM$ such that $E_2=JE_1$,
$E_4=JE_3$. According to (\ref{vtr-sympl}) (and in view of
(\ref{rho-star})) the latter identities are equivalent to
$\rho^{\ast}(X,Y)=\rho^{\ast}(Y,X)$ for every $X,Y\in T_pM$.

Now assume that $N_p\neq 0$. Then there exists $\tau\in
\Lambda^2_{+}T_pM$, $||\tau||=1$, such that $\tau\perp{\frak
J}(p)$ and $N(\tau)\neq 0$. Take a unit vector $E_1\in T_pM$ and
set $E_2=JE_1$, $E_3=K_{\tau}E_1$, $E_4=K_{{\frak
J}(p)\times\tau}E_1$. Then $E_1,...,E_4$ is an orthonormal basis
of $T_pM$ such that ${\frak J}(p)=s_1$, $\tau=s_2$, ${\frak
J}(p)\times\tau=s_3$. By (\ref{prop-N}),
$N(\tau)=N(s_2)=2N(E_1,E_3)$, $N({\frak
J}(p)\times\tau)=N(s_3)=2JN(E_1,E_3)$, thus $N({\frak
J}(p)\times\tau)=JN(\tau)$. Now we set
$A_1=||N(\tau)||^{-1}N(\tau)$, $A_2=JA_1$, $A_3=K_{\tau}A_1$,
$A_4=K_{{\frak J}(p)\times\tau}A_1$. Note that $A_4=JA_3$ by
(\ref{com}). In view of (\ref{nalpha-sympl}), we have for every
$X\in T_pM$
$$
\begin{array}{c}
g(\nabla_{X}{\frak J},\nabla_{A_1}{\frak J})\\[6pt]
=\displaystyle{\frac{1}{16}}[g(N(\tau),JX)g(N(\tau),JA)+g(N({\frak J}(p)\times\tau),JX)g(N({\frak J}(p)\times\tau),JA)]\\[6pt]
=\displaystyle{\frac{1}{16}}||N(\tau)||^2g(A_1,X).
\end{array}
$$
and
$$
g(\nabla_{X}{\frak J},\nabla_{A_2}{\frak
J})=\displaystyle{\frac{1}{16}}||N(\tau)||^2g(A_2,X).
$$
Note also that
$$
g(\nabla_{A_1}{\frak J},\nabla_{A_2}{\frak J})=0. $$ Therefore $$
(A_1)^{h}_{{\frak J}(p)}-\frac{16}{t||N(\tau)||^2}\nabla_{A_1}{\frak
J},\quad (A_2)^{h}_{{\frak
J}(p)}-\frac{16}{t||N(\tau)||^2}\nabla_{A_2}{\frak J}
$$
is a $h_t$- orthogonal basis of the normal space of ${\frak J}(M)$
at ${\frak J}(p)$. It follows from Proposition~~\ref{covder-dif}
that the normal component of the vertical part of
$Trace\,\widetilde\nabla{\frak J}_{\ast}$ at $s(p)$ vanishes if and
only if
\begin{equation}\label{vtr-l}
g(Trace\,\nabla^2{\frak J},\nabla_{A_1}{\frak
J})=g(Trace\,\nabla^2{\frak J},\nabla_{A_2}{\frak J})=0.
\end{equation}
Using (\ref{nalpha-sympl}), we see that
$$
\begin{array}{c}
\nabla_{A_1}{\frak J}=\displaystyle{\frac{1}{4}||N(\tau)||{\frak J}(p)\times\tau=\frac{1}{4}||N(\tau)||(A_1\wedge A_4+A_2\wedge A_3)},\\[6pt]
\nabla_{A_2}{\frak
J}=-\displaystyle{\frac{1}{4}||N(\tau)||\tau=-\frac{1}{4}||N(\tau)||(A_1\wedge
A_3+A_4\wedge A_2)}. \end{array}
$$
It follows from (\ref{tr-sympl}) that (\ref{vtr-l}) is equivalent to
the identities
$$
\rho^{\ast}(A_1,A_3)=\rho^{\ast}(A_2,A_4),\quad \rho^{\ast}(A_1,A_4)=\rho^{\ast}(A_4,A_1).
$$
Now, taking into account (\ref{rho-star}), we see that
(\ref{vtr-l}) is equivalent to
$\rho^{\ast}(N(\tau),X)=\rho^{\ast}(X,N(\tau))$ for
$X\perp\{N(\tau),JN(\tau)\}$, i.e.
$\rho^{\ast}(X,Y)=\rho^{\ast}(Y,X)$ for $X\perp{\cal N}_p$,
$Y\in{\cal N}_p$. The subspace ${\cal N}_p$ and ${\cal
N}_p^{\perp}$ of $T_pM$ are two-dimensional and $J$-invariant, and
it follows from (\ref{rho-star}) that
$\rho^{\ast}(X,Y)=\rho^{\ast}(Y,X)$ for $X,Y\in{\cal N}_p$ or
$X,Y\in{\cal N}_p^{\perp}$. Thus identity (\ref{vtr-l}) is
equivalent to $\rho^{\ast}$ being symmetric.

In view of  (\ref{nhtr-sympl}), the normal component of the
horizontal part of $Trace\,\widetilde\nabla{\frak J}_{\ast}$
vanishes if and only if
$$
Trace\,\{\Lambda^2_{0}TM\ni\tau\to R(\tau)(N(\tau))\}\perp\{A_1,A_2\}.
$$

This proves the statement.

\section{Examples}

In this section we give examples of almost Hermitian structures that
determine harmonic maps into twistor spaces. We also provide an
example of a Hermitian structure that is a minimal isometric
imbedding and a harmonic section of the twistor space (in sense of
C. Wood \cite{W1,W2}) but  not a harmonic map.

\subsection {Kodaira surfaces}

Recall that every primary Kodaira surface $M$ can be obtained in the
following way \cite[p.787]{Kodaira}. Let $\varphi_k(z,w)$ be the
affine transformations of ${\Bbb C}^2$ given by
$$\varphi_k(z,w) = (z+a_k,w+\overline{a}_kz+b_k),$$ where $a_k$, $b _k$, $k=1,2,3,4$, are complex
numbers such that
$$
a_1=a_2=0, \quad Im(a_3{\overline a}_4) =m b_1\neq 0,\quad b_2\neq 0
$$
for some integer $m>0$. They generate a group $G$ of transformations acting freely and properly
discontinuously on ${\Bbb C}^2$, and $M$ is the  quotient space ${\Bbb C}^2/G$.

It is well-known that $M$ can also be described  as the quotient of
${\Bbb C}^2$ endowed with a group structure by a discrete subgroup
$\Gamma$. The multiplication on ${\Bbb C}^2$ is defined by
$$
(a,b).(z,w)=(z+a,w+\overline{a}z+b),\quad (a,b), (z,w)\in  {\Bbb C}^2,
$$
and $\Gamma$ is the subgroup generated by $(a_k,b_k)$, $k=1,...,4$ (see, for example, \cite{Borc}).

Further we consider $M$ as the quotient of the group ${\Bbb C}^2$
by the discrete subgroup $\Gamma$. Every left-invariant object on
${\Bbb C}^2$ descends to a globally defined object on $M$ and both
of them will be denoted by the same symbol.

We identify ${\Bbb C}^2$ with ${\Bbb R}^4$ by $(z=x+iy,w=u+iv)\to (x,y,u,v)$ and set
$$
A_1=\frac{\partial}{\partial x}-x\frac{\partial}{\partial
u}+y\frac{\partial}{\partial v},\quad A_2=\frac{\partial}{\partial
y}-y\frac{\partial}{\partial u}-x\frac{\partial}{\partial v},\quad
A_3=\frac{\partial}{\partial u},\quad A_4=\frac{\partial}{\partial
v}.
$$
These form a basis for the space of left-invariant vector fields
on ${\Bbb C}^2$. We note that their Lie brackets are
$$
[A_1,A_2]=-2A_4,\quad [A_i,A_j]=0
$$
for all other $i,j$. It follows that the group ${\Bbb C}^2$ defined
above is solvable.

 Denote by $g$ the left-invariant Riemannian metric on $M$ for which the basis
$A_1,...,A_4$ is orthonormal.

We shall show that any integrable or symplectic almost complex
structures $J$ on $M$ compatible with the metric $g$ and defined by
a left-invariant almost complex structure on ${\Bbb C}^2$ is a
harmonic map from $(M,g)$ to $({\cal Z},h_t)$.

Note that by \cite{Has} every complex structure on $M$ is induced by
a left-invariant complex structure on ${\Bbb C}^2$.

\medskip

\noindent {\bf I}. If $J$ is a left-invariant almost complex
structure compatible with $g$, we have $JA_i=\sum_{j=1}^4 a_{ij}A_j$
where $a_{ij}$ are constants with $a_{ij}=-a_{ji}$.  Let $N$ be the
Nijenhuis tensor of $J$. Computing $N(A_i,A_j)$ in terms of
$a_{ij}$, one can see (\cite{M,D}) that $J$ is integrable if and
only if
$$
JA_1=\varepsilon_1 A_2,\quad JA_3=\varepsilon_2 A_4,\quad \varepsilon_1,\varepsilon_2=\pm 1.
$$
Denote by $\theta$ the Lee form of the Hermitian structure
$(M,g,J)$ where $J$ is defined by means of the latter identities.

The non-zero covariant derivatives $\nabla_{A_i}A_j$ are
$$
\nabla_{A_1}A_2=-\nabla_{A_2}A_1=-A_4,\quad \nabla_{A_1}A_4=\nabla_{A_4}A_1=A_2,\quad
\nabla_{A_2}A_4=\nabla_{A_4}A_2=-A_1.
$$
This implies that the Lie form is
$$
\theta(X)=-2\varepsilon_1g(X,A_3).
$$
Therefore
\begin{equation}\label{cond-1}
B=-2\varepsilon_1A_3,\quad \nabla\theta=0.
\end{equation}
Set for short $R_{ijk}=R(A_1,A_j)A_k$. Then the non-zero $R_{ijk}$
are
$$
\begin{array}{c}
R_{121}=-3A_2,\quad R_{122}=3A_1,\quad R_{141}=A_4,\\[6pt]
 R_{144}=-A_1,\quad R_{242}=A_4,\quad R_{244}=-A_2.
\end{array}
$$
Set also $\rho_{ij}=\rho(A_i,A_j)$, $\rho^{\ast}_{ij}=\rho^{\ast}(A_i,A_j)$. Then
\begin{equation}\label{cond-2}
\begin{array}{c}
\rho_{ij}=0\quad {\mbox and}\quad \rho^{\ast}_{ij}=0\quad {\mbox except}\\[6pt]
\rho_{11}=\rho_{22}=-2,\quad \rho_{44}=2,\quad \rho^{\ast}_{11}=\rho^{\ast}_{22}=-3.
\end{array}
\end{equation}
It follows form (\ref{cond-1}), (\ref{cond-2}) and
Theorem~\ref{harm-int} that the complex structure $J$ is a harmonic
map from $(M,g)$ to the twistor space $({\cal Z},h_t)$.

It is easy to describe explicitly the twistor space  $({\cal
Z},h_t)$  (\cite{D}) since $\Lambda^2_{+}M$ admits a global
orthonormal frame defined by
$$
s_1=\varepsilon_1 A_1\wedge A_2+\varepsilon_2A_3\wedge A_4,\quad s_2=A_1\wedge A_3+\varepsilon_1\varepsilon_2 A_4\wedge A_2,\quad s_3=\varepsilon_2 A_1\wedge A_4+\varepsilon_1 A_2\wedge A_3.
$$
It induces a natural diffeomorphism $F: {\cal Z}\cong M\times S^2$,
$\sum_{k=1}^3 x_k s_k(p)\to (p,x_1,x_2,x_3)$, under which $J$
determines the section $p\to (p,1,0,0)$. In order to find an
explicit formula for the metrics $h_t$ we need the covariant
derivatives of $s_1,s_2,s_3$ with respect to the Levi-Civita
connection $\nabla$ of $g$. The non-zero of these are
$$
\begin{array}{c}
\nabla_{A_1}s_1=-\epsilon_2\nabla_{A_4}s_2=-\varepsilon_1\varepsilon_2 s_3,\quad
\varepsilon_1 \nabla_{A_1}s_3=-\nabla_{A_2}s_2=\varepsilon_2s_1,\\[6pt]
\varepsilon_2\nabla_{A_2}s_1=-\varepsilon_1\nabla_{A_4}s_3=s_2.
\end{array}
$$
It follows that $F_{\ast}$ sends the horizontal lifts
$A_1^h,...,A_4^h$ at a point $\sigma=\sum_{k=1}^3 x_k
s_k(p)\in{\cal Z}$ to the following vectors of $TM\oplus TS^2$
$$
A_1+\varepsilon_1\varepsilon_2 (-x_3,0,x_1),\quad A_2+\varepsilon_2(x_2,-x_1,0),\quad A_3,\quad A_4+\varepsilon_1 (0,x_3,-x_2).
$$
For $x=(x_1,x_2,x_3)\in S^2$, set
$$
u_1(x)=\varepsilon_1\varepsilon_2 (-x_3,0,x_1),\quad u_2(x)=\varepsilon_2(x_2,-x_1,0),\quad
u_3(x)=0,\quad u_4(x)=\varepsilon_1 (0,x_3,-x_2).
$$
Denote the pushforward of the metric $h_t$ under $F$ again by
$h_t$. Then, if $X,Y\in T_pM$ and $P,Q\in T_xS^2$,
\begin{equation}\label{ht}
h_t(X+P,Y+Q)=g(X,Y)+t<P-\sum_{i=1}^4 g(X,A_i)u_i(x),Q-\sum_{j=1}^4
g(Y,A_j)u_j(x)>
\end{equation}
where $<.,.>$ is the standard metric of ${\Bbb R}^3$.

\medskip

\noindent {\bf II}. Suppose again that $J$ is an almost complex
structure on $M$ obtained from a left-invariant almost complex
structure on ${\Bbb C}^2$ and compatible with the metric $g$. Set
$JA_i=\sum_{j=1}^4 a_{ij}A_j$. Denote  the fundamental $2$-form of
the almost Hermitian structure $(g,J)$ by $\Omega$.  The basis
dual to $A_1,...,A_4$ is
$$
{\alpha}_1=dx,\quad {\alpha}_2=dy,\quad
{\alpha}_3=xdx+ydy+du,\quad {\alpha}_4=-ydx+xdy+dv.
$$
 We have $d{\alpha}_1=d{\alpha}_2=d{\alpha}_3=0,\> d{\alpha}_4=2dx\wedge
 dy$.
Hence
$$
d\Omega=d\sum_{i<j}a_{ij}{\alpha}_i\wedge{\alpha}_j=-2a_{34}dx\wedge
dy\wedge du.$$ Thus $d\Omega=0$ is equivalent to $a_{34}=0$. If
$a_{34}=0$, we have $a_{1j}=a_{2j}=0$ for $j=1,2$,
$a_{3k}=a_{4k}=0$ for $k=3,4$, $a_{13}^2+a_{14}^2=1$,
$a_{13}a_{23}+a_{14}a_{24}=0$, $a_{23}^2+a_{24}^2=1$. It follows
that  the structure $(g,J)$ is almost K\"ahler (symplectic) if and
only if $J$ is given by (\cite{M,D})
$$
\begin{array}{c}
JA_1=-\varepsilon_1\sin\varphi A_3+\varepsilon_1\varepsilon_2\cos\varphi A_4,\quad
JA_2=-\cos\varphi A_3-\varepsilon_2\sin\varphi A_4, \\[6pt]
JA_3=\varepsilon_1\sin\varphi A_1+\cos\varphi A_2,\quad JA_4=-\varepsilon_1\varepsilon_2\cos\varphi
A_1+\varepsilon_2\sin\varphi A_2, \\[6pt]
\varepsilon_1,\varepsilon_2=\pm 1, \> \varphi\in [0,2\pi).
\end{array}
$$

Let $J$ the almost complex structure defined by these identities
for some $\varepsilon_1,\varepsilon_2,\varphi$. Set
$$
E_1=A_1,\quad E_2=-\varepsilon_1\sin\varphi A_3+\varepsilon_1\varepsilon_2\cos\varphi A_4,\quad
E_3=\cos\varphi A_3+\varepsilon_2\sin\varphi A_4,\quad E_4=A_2.
$$
Then $E_1,...,E_4$ is an orthonormal frame of $TM$ for which
$JE_1=E_2$ and $JE_3=E_4$. The only non-zero Lie bracket of these
fields is
$$
[E_1,E_4]=-2(\varepsilon_1\varepsilon_2\cos\varphi E_2+\varepsilon_2\sin\varphi E_3).
$$
The non-zero covariant derivatives $\nabla_{E_i}E_j$ are
$$
\begin{array}{c}
\nabla_{E_1}E_2=\nabla_{E_2}E_1=\varepsilon_1\varepsilon_2\cos\varphi E_4,\quad
\nabla_{E_1}E_3=\nabla_{E_3}E_1=\varepsilon_2\sin\varphi E_4,\\[6pt]
\nabla_{E_1}E_4=-\nabla_{E_4}E_1=-\varepsilon_1\varepsilon_2\cos\varphi
E_2-\varepsilon_2\sin\varphi
E_3,\\[6pt]
\nabla_{E_2}E_4=\nabla_{E_4}E_2=-\varepsilon_1\varepsilon_2\cos\varphi E_1,\quad
\nabla_{E_3}E_4=\nabla_{E_4}E_3=-\varepsilon_2\sin\varphi E_1.
\end{array}
$$
Set $R_{ijk}=R(E_i,E_j)E_k$. We have the following table for the
non-zero components of the curvature tensor $R$:
$$
\begin{array}{c}
R_{121}=\cos^2\varphi E_2+\displaystyle{\frac{1}{2}}\varepsilon_1\sin 2\varphi E_3,\quad R_{122}=-\cos^2\varphi E_1,\quad
R_{123}=-\displaystyle{\frac{1}{2}}\varepsilon_1\sin 2\varphi E_1, \\[6pt]
R_{131}=\displaystyle{\frac{1}{2}}\varepsilon_1\sin 2\varphi E_2+\sin^2\varphi E_3,\quad R_{132}=-\displaystyle{\frac{1}{2}}\varepsilon_1\sin 2\varphi E_1,
\quad R_{133}=-\sin^2\varphi E_1,\\[6pt]
R_{141}=-3E_4,\quad R_{144}=3E_1,\\[6pt]
R_{242}=\cos^2\varphi E_4,\quad R_{243}=\displaystyle{\frac{1}{2}}\varepsilon_1\sin 2\varphi E_4,
\quad R_{244}=-\cos^2\varphi E_2-\displaystyle{\frac{1}{2}}\varepsilon_1\sin 2\varphi E_3,\\[6pt]
R_{342}=\displaystyle{\frac{1}{2}}\varepsilon_1\sin 2\varphi
E_4,\quad R_{343}=\sin^2\varphi E_4, \quad
R_{344}=-\displaystyle{\frac{1}{2}}\varepsilon_1\sin 2\varphi
E_2-\sin^2\varphi E_3.
\end{array}
$$
Define an orthonormal frame $s_l$, $l=1,2,3$, of $\Lambda^2_{+}TM$
by means of $E_1,...,E_4$. Then  $s_2,s_3$ is a frame of
$\Lambda^2_{0}TM$ and by (\ref{prop-N})
$$
\begin{array}{c}
N(s_2)=2N(E_1,E_3)=-4\varepsilon_1\varepsilon_2\cos\varphi E_1+4\varepsilon_2\sin\varphi E_4,\\[6pt]
N(s_3)=2N(E_1,E_4)= 4\varepsilon_1\varepsilon_2\cos\varphi E_2+4\varepsilon_2\sin\varphi E_3.
\end{array}
$$
It follows that
$$
Trace\,\{\Lambda^2_{0}TM\ni\tau\to
R(\tau)(N(\tau))\}=R(s_2)(N(s_2))+R(s_3)(N(s_3))=0.
$$

Setting $\rho^{\ast}_{ij}=\rho^{\ast}(E_i,E_j)$, we have
$$
\begin{array}{c}
\rho^{\ast}_{11}=\rho^{\ast}_{22}=\cos^2\varphi,\quad
\rho^{\ast}_{33}=\rho^{\ast}_{44}=\sin^2\varphi ,\quad
\rho^{\ast}_{14}=\rho^{\ast}_{41}= -\displaystyle{\frac{1}{2}}\varepsilon_1\sin 2\varphi,\\[6pt]
\end{array}
$$
and the other $\rho^{\ast}_{ij}$ vanish.  Thus, by
Theorem~\ref{harm-sympl}, the almost K\"ahler structure $J$ is a
harmonic map $(M,g)\to ({\cal Z},h_t)$.

As in the preceding case, it is easy to find an explicit description
of the twistor space ${\cal Z}$ of $M$ and the metric $h_t$
(\cite{D}). The frame $s_1,s_2,s_3$ gives rise to an obvious
diffeomorphism $F: {\cal Z}\cong M\times S^2$ under which ${\frak
J}$ becomes the map $p\to (p,1,0,0)$. We have the following table
for the covariant derivatives of $s_1,s_2,s_3$:
$$
\begin{array}{c}
\nabla_{E_1}s_1=\nabla_{E_4}s_2=\varepsilon_1\varepsilon_2\cos\varphi s_3,\quad
\nabla_{E_4}s_1=-\nabla_{E_1}s_2=-\varepsilon_2\sin\varphi s_3,\\[6pt]
\nabla_{E_2}s_1=\varepsilon_1\varepsilon_2\cos\varphi s_2,\quad
\nabla_{E_2}s_2=-\varepsilon_1\varepsilon_2\cos\varphi s_1,\quad
\nabla_{E_2}s_3=0,\\[6pt]
\nabla_{E_3}s_1=\varepsilon_2\sin\varphi s_2,\quad
\nabla_{E_3}s_2=-\varepsilon_2\sin\varphi s_1,\quad \nabla_{E_3}s_3=0,\\[6pt]
\nabla_{E_1}s_3=-\varepsilon_1\varepsilon_2\cos\varphi s_1-\varepsilon_2\sin\varphi
s_2,\quad \nabla_{E_4}s_3=\varepsilon_2\sin\varphi s_1
-\varepsilon_1\varepsilon_2\cos\varphi s_2.
\end{array}
$$
Using this table we see that $F_{\ast}$ sends the horizontal lifts
$E_i^h$, $i=1,...,4$, to $E_i+u_i$ where
$$
\begin{array}{c}
u_1(x)=(x_3\varepsilon_1\varepsilon_2\cos\varphi, x_3\varepsilon_2\sin\varphi,-x_1\varepsilon_1\varepsilon_2\cos\varphi-x_2\varepsilon_2\sin\varphi),\\[6pt]
u_2(x)=(x_2\varepsilon_1\varepsilon_2\cos\varphi,
-x_1\varepsilon_1\varepsilon_2\cos\varphi, 0) ,\quad
u_3(x)=(x_2\varepsilon_2\sin\varphi,-x_1\varepsilon_2\sin\varphi,0),\\[6pt]
u_4(x)=(-x_3\varepsilon_2\sin\varphi,
x_3\varepsilon_1\varepsilon_2\cos\varphi,x_1\varepsilon_2\sin\varphi-x_2\varepsilon_1\varepsilon_2\cos\varphi).
\end{array}
$$
for $x=(x_1,x_2,x_3)\in S^2$. Then, if $X,Y\in T_pM$ and $P,Q\in T_xS^2$,
\begin{equation}\label{ht-sympl}
h_t(X+P,Y+Q)=g(X,Y)+t<P-\sum_{i=1}^4 g(X,E_i)u_i(x),Q-\sum_{j=1}^4
g(Y,E_j)u_j(x)>.
\end{equation}

\subsection{Four-dimensional Lie groups} We shall show that every
left-invariant almost K\"ahler structure $(g,J)$ with
$J$-invariant Ricci tensor on a $4$-dimensional Lie group $M$
determines a harmonic map ${\frak J}: (M,g)\to({\cal Z},h_t)$.

These (non-integrable) structures have been determined in \cite{F}.
According to the main result therein, for any such a structure
$(g,J)$, there exists an orthonormal frame of left-invariant vector
fields $E_1,...,E_4$ such that
$$
JE_1=E_2,\quad JE_3=E_4
$$
and
$$
\begin{array}{c}
\displaystyle{[E_1,E_2]=0,\quad [E_1,E_3]=sE_1+\frac{s^2}{t}E_2,\quad [E_1,E_4]=\frac{s^2-t^2}{2t}E_1-sE_2},\\[6pt]
\displaystyle{[E_2,E_3]=-tE_1-sE_2,\quad [E_2,E_4]=-sE_1-\frac{s^2-t^2}{2t}E_2,\quad [E_3,E_4]=-\frac{s^2+t^2}{t}E_3},
\end{array}
$$
where $s$ and $t\neq 0$ are real numbers. Then we have the following table for the Levi-Civita connection
$$
\begin{array}{c}
\displaystyle{\nabla_{E_1}E_1=-sE_3-\frac{s^2-t^2}{2t}E_4, \quad \nabla_{E_2}E_1=-\frac{s^2-t^2}{2t}E_3+sE_4, \quad \nabla_{E_3}E_1=\frac{s^2+t^2}{2t}E_2},\\[8pt]
\displaystyle{\nabla_{E_1}E_2=-\frac{s^2-t^2}{2t}E_3+sE_4,~\quad \nabla_{E_2}E_2=sE_3+\frac{s^2-t^2}{2t}E_4,~\quad \nabla_{E_3}E_2==\frac{s^2+t^2}{2t}E_1},\\[8pt]
\displaystyle{\nabla_{E_1}E_3=sE_1+\frac{s^2-t^2}{2t}E_2,~\quad \nabla_{E_2}E_3=\frac{s^2-t^2}{2t}E_1-sE_2,~\quad\nabla_{E_3}E_3=\frac{s^2+t^2}{t}E_4},\\[8pt]
\displaystyle{\nabla_{E_1}E_4=\frac{s^2-t^2}{2t}E_1-sE_2,~\quad \nabla_{E_2}E_4=-sE_1-\frac{s^2-t^2}{2t}E_2,~\quad\nabla_{E_3}E_4=-\frac{s^2+t^2}{t}E_3},\\[8pt]
\nabla_{E_4}E_1=\nabla_{E_4}E_2=\nabla_{E_4}E_3=\nabla_{E_4}E_4=0.
\end{array}
$$
This implies the following table for the components
$R_{ijk}=R(E_i,E_j)E_k$ of the  curvature tensor; in this table
$\lambda=\displaystyle{\frac{s^2+t^2}{2t}}$. $$
\begin{array}{c}
R_{121}=2\lambda E_2, \quad R_{122}=-2\lambda E_1, \quad R_{123}=2\lambda E_4,
\quad R_{124}=-2\lambda E_3,\\[6pt]
R_{131}=-\lambda E_3, \quad R_{132}=\lambda E_4, \quad R_{133}=\lambda E_1, \quad R_{134} =-\lambda E_2, \\[6pt]
R_{141}=-\lambda E_4, \quad R_{142}=-\lambda E_3, \quad R_{143}=\lambda E_2, \quad R_{144}=\lambda E_1, \\[6pt]
R_{231}=-\lambda E_4, \quad R_{232}=-\lambda E_3, \quad R_{233}=\lambda E_2, \quad R_{234}=\lambda E_1, \\[6pt]
R_{241}=\lambda E_3, \quad R_{242}=-\lambda E_4, \quad R_{243}=-\lambda E_1, \quad R_{244}=\lambda E_2, \\[6pt]
R_{341}=2\lambda E_2, \quad R_{342}=-2\lambda E_1, \quad R_{343}=-4\lambda E_4, \quad R_{344}=4\lambda E_3.
\end{array}
$$
Then the non-zero $\rho^{\ast}_{ij}=\rho^{\ast}(E_i,E_j)$ are
$$
\rho^{\ast}_{11}=\rho^{\ast}_{22}=4\lambda,\quad
\rho^{\ast}_{33}=\rho^{\ast}_{44}=-2\lambda.
$$
Therefore the $\ast$-Ricci tensor $\rho^{\ast}$ is symmetric.

Set $s_2=E_1\wedge E_3+E_4\wedge E_2$, $s_3=E_1\wedge E_4+E_2\wedge
E_3$. Then
$$
\begin{array}{c}
N(s_2)=2N(E_1,E_3)=-8(sE_1+\displaystyle{\frac{s^2-t^2}{2t}}E_2), \\[6pt] N(s_3)=2N(E_1,E_4)
=8(-\displaystyle{\frac{s^2-t^2}{2t}}E_1+sE_2).
\end{array}
$$
It follows that
$$
Trace\,\{\Lambda^2_{0}TM\ni\tau\to R(\tau)(N(\tau))\}=0.
$$
Thus, according to Theorem~\ref{min-sympl}, $J$ is a harmonic map.

\subsection{ Inoue surfaces of type $S^0$} It has been observed in
\cite{Tr} that every Inoue surface $M$ of type $S^0$ admits a
locally conformal K\"ahler metric $g$ (cf. also \cite{DO}). The map
${\frak J}: (M,g)\to({\cal Z},h_t)$ determined by the complex
structure of $M$ is a minimal isometric imbedding that is a harmonic
section of the twistor space but  not a harmonic map. To see this
let us first recall the construction of the Inoue surfaces of type
$S^0$ (\cite{Inoue}). Let $A\in SL(3,\mathbb{Z})$ be a matrix with a
real eigenvalue $\alpha > 1$ and two complex eigenvalues $\beta$ and
$\overline{\beta}$, $\beta\neq\overline{\beta}$. Choose eigenvectors
$(a_1,a_2,a_3)\in{\mathbb R}^3$ and $(b_1,b_2,b_3)\in {\mathbb C}^3$
of $A$ corresponding to $\alpha$ and $\beta$, respectively. Then the
vectors $(a_1,a_2,a_3), (b_1,b_2,b_3),
(\overline{b_1},\overline{b_2},\overline{b_3})$  are
$\mathbb{C}$-linearly independent. Denote the upper-half plane in
$\mathbb{C}$ by ${\bf H}$ and let $\Gamma$ be the group of
holomorphic automorphisms of ${\bf H}\times {\mathbb C}$ generated
by
$$g_o:(w,z)\to (\alpha w,\beta z), \quad g_i:(w,z)\to (w+a_i,z+b_i), \>i=1,2,3 .$$
The group $\Gamma$ acts on ${\bf H}\times{\mathbb C}$ freely and
properly discontinuously.  Then $M=({\bf H}\times {\mathbb
C})/\Gamma$ is a complex surface known as Inoue surface of type
$S^0$. As in \cite{Tr}, consider on ${\bf H}\times {\mathbb C}$ the
Hermitian metric
$$
g=\frac{1}{v^2}(du\otimes du+dv\otimes dv)+v(dx\otimes dx+dy\otimes
dy),\quad u+iv\in{\bf H}, \quad x+iy\in {\mathbb C}.
$$
This metric is invariant under the action of the group $\Gamma$, so
it descends to a Hermitian metric on $M$ which we denote again by
$g$. Instead on $M$, we shall work with $\Gamma$-invariant objects
on ${\bf H}\times{\mathbb C}$.  Let $\Omega$ be the fundamental
$2$-form of the Hermitian structure $(g,J)$ on ${\bf
H}\times{\mathbb C}$, $J$ being the standard complex structure. Then
$$
d\Omega=\frac{1}{v}dv\wedge\Omega.
$$
Hence the Lee form is $\theta=d\ln v$. In particular, $d\theta=0$,
i.e. $(g,J)$ is a locally conformal K\"ahler structure. Set
$$
E_1=v\frac{\partial}{\partial u},\quad E_2=v\frac{\partial}{\partial
v},\quad E_3=\frac{1}{\sqrt v}\frac{\partial}{\partial x},\quad
E_4=\frac{1}{\sqrt v}\frac{\partial}{\partial y}.
$$
These are $\Gamma$-invariant vector fields constituting an
orthonormal basis. Note that the vector field dual to the Lee form
is $B=E_2$. The non-zero Lie brackets of $E_1,...,E_4$ are
$$
[E_1,E_2]=-E_1,\quad [E_2,E_3]=-\frac{1}{2}E_3,\quad
[E_2,E_4]=-\frac{1}{2}E_4.
$$
Then we have the following table for the Levi-Civita connection
$\nabla$ of $g$:
$$
\begin{array}{c}
\nabla_{E_1}E_1=E_2,\quad \nabla_{E_1}E_2=-E_1,\\[6pt]
\displaystyle{\nabla_{E_3}E_2=\frac{1}{2}E_3,\quad
\nabla_{E_3}E_3=-\frac{1}{2}E_2, \quad
\nabla_{E_4}E_2=\frac{1}{2}E_4,\quad
\nabla_{E_4}E_4=-\frac{1}{2}E_2},
\end{array}
$$
and all other $\nabla_{E_i}E_j=0$. It follows that
$$
\rho(E_k,B)=\rho^{\ast}(E_k,B)=0~~\rm{for}~k=1,3,4,\quad
\rho(E_2,B)=\frac{1}{2},\quad \rho^{\ast}(E_2,B)=-1.
$$
By Corollary~\ref{Wood} ${\frak J}:(M,g)\to ({\cal Z},h_t)$ is a
harmonic section. It is also a minimal isometric imbedding by
Theorem~\ref{min-int}. But ${\frak J}$ is not a harmonic map
according to Theorem~\ref{harm-int}.

\end{document}